\documentclass{amsart}
\usepackage{epsf,amssymb,xypic,amsthm}
\begin{document}

\newcommand{\mT}{\mathcal{T}}
\newcommand{\mB}{\mathcal{B}}
\renewcommand{\Delta}{\Phi}
\newcommand{\zeroh}{\hat{0}}
\newcommand{\oneh}{\hat{1}}
\newcommand{\piu}{p^\uparrow}
\newcommand{\pid}{p^\downarrow}
\newcommand{\Id}{I^\downarrow}
\newcommand{\Iu}{I^\uparrow}

\newtheorem{theorem}{Theorem}
\newtheorem{proposition}{Proposition}
\newtheorem{lemma}{Lemma}
\newtheorem{conjecture}{Conjecture}
\theoremstyle{remark}
\newtheorem*{remark}{Remark}
\newtheorem*{example}{Example}
\theoremstyle{definition}
\newtheorem*{definition}{Definition}

\title [Analogue of distributivity]
{An analogue of distributivity for ungraded lattices}

\author {Hugh Thomas}

\date {September 2005}

\address 
{Department of Mathematics and Statistics, University of New Brunswick,
Fredericton NB, E3B 5A3 Canada}

\keywords{ 
Left modular lattices, supersolvable lattices, 
Tamari lattice, Cambrian lattice}

\email{hugh@math.unb.ca}
\begin{abstract}
In this paper, we define a property, trimness, for lattices.  Trimness is 
a not-necessarily-graded generalization of distributivity; in particular,
if a lattice is trim and graded, it is distributive.  
Trimness is preserved under taking 
intervals and suitable sublattices.  
Trim lattices satisfy a weakened form of modularity.  
The order complex of a trim lattice
is contractible or homotopic to a sphere; the latter holds exactly if the
maximum element of the lattice is a join of atoms.

Other than distributive
lattices, the main examples of trim lattices are the Tamari lattices
and various generalizations of them.  We show that the Cambrian lattices 
in types $A$ and $B$ defined by Reading 
are trim, and we conjecture that all Cambrian lattices are trim.  

\end{abstract}
\maketitle

\section{Introduction}

Some of the first examples of lattices which anyone encounters are the
finite distributive lattices.  Supersolvable lattices are a generalization
of them introduced by Stanley [St] in 1972.  
Lattices of both these types are necessarily graded.  
Left modular lattices were introduced by Blass and Sagan [BS] as a further
generalization of supersolvable lattices.  
In [MT], combining results from [Mc] and [Li], McNamara and the author 
showed that left modularity for lattices can be thought of as 
``supersolvability without gradedness,''  
in the sense that supersolvable lattices are
left modular (as was shown in [St]), and all graded left modular 
lattices are supersolvable [MT].  
Thus, we have the following diagram:

$$
\epsfbox{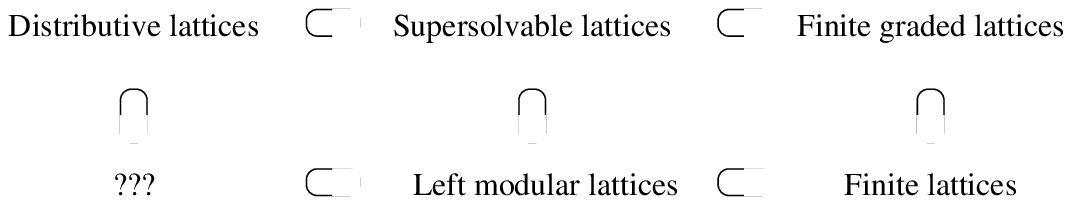}
$$
$$\text{Figure 1}$$

This paper is an attempt to provide something to fit in place of the ???
in the diagram, that is to say, a not-necessarily-graded generalization
of distributivity.  

We begin with some necessary definitions.  All our lattices are assumed
to be finite.  
An element $x$ of a lattice $L$ is said to be {\it left modular} if for any 
$y<z$ in $L$, 
$$(y\vee x)\wedge z=y\vee(x\wedge z).$$
A set of pairwise comparable elements is called a {\it chain}.  A {\it maximal
chain} (also called an unrefinable chain) is one which is maximal with
respect to inclusion.  A lattice is called {\it graded} if every maximal
chain in a given interval is of the same cardinality.  
A lattice is said to be {\it left modular} if it has a 
maximal chain of elements all of which are left modular.  

An element in a lattice is called {\it join-irreducible} if it cannot be 
written as the join of two strictly smaller elements.  ($\zeroh$, the 
minimum element of a lattice, does not count as a join-irreducible.)
Dually, an
element is called {\it meet-irreducible} if it cannot be written
as the meet of two strictly larger elements.  
(Similarly, $\oneh$, the maximum element of a lattice, does not count
as a meet-irreducible.) 
If a lattice has a maximal chain of $n+1$ elements, then it must have at 
least $n$ join-irreducibles and at least $n$ meet-irreducibles.    

\begin{definition}
We say that a lattice is {\it trim} if it has a maximal chain of
$n+1$ left modular elements, exactly $n$ join-irreducibles, and
exactly $n$ meet-irreducibles (that is to say, the minimum possible
number of each).  
\end{definition}

Distributive lattices are example of trim lattices, though not
typical examples, since distributive lattices are graded and trim lattices
need not be.  However, if a lattice is graded and trim, then it is 
distributive (Theorem 2).  This is a special case of a theorem of Markowsky
[Ma].    

In this paper, we investigate some of the properties of trim lattices,
and show that they are in many respects similar to distributive
lattices.  
We show that if $L$ is trim, then so are its intervals (Theorem 1), and
so are its sublattices which contain its left modular chain
(Theorem 3).  We also show that if $G$ is a group which acts on a trim
lattice $L$
by lattice automorphisms, then $L^G$, the sublattice of $L$ consisting
of elements fixed by
$G$, is again trim (Theorem 4).

We show that a trim lattice satisfies the level condition
of [BS] (Theorem 5).  
Left modular lattices
satisfying this condition are known as LL-lattices [BS]; thus, trim
lattices are LL-lattices.  

One consequence of 
Theorem 5 is that in a trim lattice,
if $w$ is covered by $y$ and $z$, then $y\vee z$ covers at least one of
$y$ and $z$  (Theorem 6).  This is a weakened form of modularity.  
Another consequence of Theorem 5 is 
that the order complex of a trim lattice is either
contractible or homotopic to a sphere, and the latter holds exactly
when the maximum element of the lattice is a join of atoms (Theorem 7).  
(Note that since intervals of trim
lattices are trim, Theorem 7 can also be applied to any interval in a 
trim lattice, thus showing that the order complex of any interval is
again either homotopic to a sphere or contractible.)

In [Re], Reading introduced a family of {\it Cambrian lattices} for each
finite reflection group.  (Recall that finite reflection groups consist of
four infinite families, $A_n$ ($n\geq 1$), $B_n$ ($n\geq 2$), 
$D_n$ ($n\geq 4$), $I_2(n)$ ($n=5$ or $n\geq 7$), and seven 
exceptional groups, $E_6$, $E_7$, $E_8$, $F_4$, $G_2$, $H_3$, $H_4$.)
The Cambrian lattices in type $A$ include the 
classical Tamari lattice (which goes back to [Ta]; a more recent reference is
[BW]) 
and in type $B$ 
include
the type $B$
Tamari lattice (also studied in [Th]).  We show that all the Cambrian lattices
in types $A$ and $B$ are trim (Theorems 8 and 9).
The Cambrian lattices in the other types are
not yet well understood, but we offer the following conjecture:

\begin{conjecture}All Cambrian lattices are trim.\end{conjecture}

\section{Left Modular Lattices}

The proofs in this paper depend on the theory of left modular 
lattices.  The study of such lattices was 
initiated by Blass and Sagan [BS], and continued in Liu [Li], 
Liu and Sagan [LS], and McNamara and Thomas [MT].  
We will begin with a review of the properties of left modular lattices.  
More details on all of these properties can be found in [MT].

\begin{proposition}[{[MT]}] If $L$ has a left modular maximal chain
$\zeroh=x_0\lessdot x_1\lessdot\dots\lessdot \oneh$, then any interval
$[y,z]$ also has a left modular maximal chain.   
More precisely, the elements $y\vee x_i\wedge z$ form a left modular
maximal chain in $[y,z]$.  (Note that the $y\vee x_i \wedge z$ will not
all be distinct.)
\end{proposition}

We define three edge-labellings of a left modular
lattice $L$ with left modular maximal 
chain $\zeroh=x_0\lessdot x_1\lessdot\dots\lessdot x_n=\oneh$, which we refer
to as the labelling induced from join-irreducibles, that induced from
meet-irreducibles, and that induced from the left modular chain.  
(Note that these labellings all depend on the prior choice of a left modular
maximal chain.)

If $v$ is 
a join-irreducible of $L$, we label it by the natural number
$$\delta(v)=\min(\{i\mid v\leq x_i\}).$$  
Now, for any $y \lessdot z$, define
$$\gamma_1(y\lessdot z)=\min(\{\delta(v)\mid v \text{ join-irreducible, }
v\leq z,\, v \not \leq y\}).$$

The labelling induced from meet-irreducibles is defined similarly.  
If $v$ is a 
meet-irreducible, we set 
$$\epsilon(v)=\max(\{i\mid x_i\leq v\})+1.$$
(Except for the $+1$, this is just the dual of the definition of 
$\delta$.)  Now, as we did for the 
labelling induced from join-irreducibles, for $y \lessdot z$, we define 
$$\gamma_2(y\lessdot z)=\max(\{\epsilon(v)\mid v\text{ meet-irreducible, }
v\geq y,\, v\not \geq z\}).$$

Thirdly, the labelling induced from the left modular chain is defined
as follows: 
$$\gamma_3(y\lessdot z)=\min(\{i\mid y\vee x_i \wedge z=z\}).$$

\begin{proposition}[{[Li]}] For any left modular lattice with a specified
left modular maximal chain, 
the three labellings $\gamma_1$, $\gamma_2$, and $\gamma_3$ coincide.
\end{proposition}

The fact that $\gamma_1$ and $\gamma_3$ coincide is proved in [Li]; 
the dual of that result shows that $\gamma_2$ and $\gamma_3$ coincide.  
Since the three labellings coincide, we will drop the subscripts and
denote the labelling by $\gamma$.

A labelling of the edges of the Hasse diagram of a poset is called an 
EL-labelling [B1] if it satisfies the following two properties:

(i) In any interval, there is a unique maximal chain which has the 
property that the labels on the chain strictly increase as you read up the
chain. (This chain is called the ``increasing chain''.)

(ii) In any interval, 
the label word obtained by reading up the increasing chain 
lexicographically precedes the word obtained by reading up any other
maximal chain in the interval.

(In our context, the labellings of the edges of a Hasse diagram will always
be positive integers with the usual order.  In general, the labels may be
drawn from any poset; this introduces some additional technicalities which
we shall not need to refer to.)

If a partially ordered set admits an EL-labelling then its order complex
is shellable, and is therefore homotopic to a wedge of spheres, one
for each maximal chain such that the labels weakly decrease as you
read up the chain.  (Such chains are called ``decreasing chains''.)
The dimension of the sphere corresponding to a given decreasing chain is
two less than the length of the chain.

\begin{proposition}[{[Li]}] For a left modular lattice $L$, 
the edge-labelling of $L$ already described 
is an EL-labelling.  \end{proposition}

In fact, we can say more about the labelling of a left modular lattice.  
In [MT], we defined {\it interpolating labellings} to be EL-labellings
such that in addition, if $v\lessdot u \lessdot w$ is a maximal chain which
is not increasing, and the corresponding increasing chain is $v=y_0\lessdot y_1
\lessdot \dots \lessdot y_r=w$, then the label of $v\lessdot u$ coincides 
with the label of $y_{r-1}\lessdot y_r$, and the label of $u\lessdot w$
coincides with the label of $y_0\lessdot y_1$.  We showed the following
proposition:

\begin{proposition}[{[MT]}] If $L$ is a left modular lattice, then the
labelling defined above is interpolating.  Conversely, if a lattice $L$
admits an interpolating labelling, then the elements of the increasing
chain from $\zeroh$ to $\oneh$ are left modular, and therefore $L$ is left 
modular.  
\end{proposition}

We need one more result from [MT] about labellings of intervals.  
Let $[y,z]$ be an interval in a left modular lattice $L$.  Since the 
$y\vee x_i \wedge z$ form a left modular chain in $[y,z]$, 
the above construction can be applied to yield an EL-labelling.
The restriction of the labelling of $L$ to $[y,z]$ also yields an EL-labelling.
These two labellings typically do not coincide for the trivial reason
that their label sets differ.  However, we have the following proposition:

\begin{proposition}[{[MT]}]
Let $[y,z]$ be an interval in a lattice $L$ with left modular maximal
chain $\zeroh=x_0\lessdot x_1\lessdot\dots\lessdot x_n=\oneh$.  The labelling of $L$ restricted to $[y,z]$ agrees 
(up to an order-preserving relabelling) with the 
labelling which $[y,z]$ has as a lattice with left modular chain
$y\vee x_i \wedge z$.
\end{proposition}

(When we speak of an order-preserving relabelling, we mean that one label
set has been replaced by a different label set, but the relative orders
of the labels have been preserved.)

We record here one additional lemma about left modular lattices which
we shall need.  

\begin{lemma} Let $L$ be a lattice
with left modular maximal chain $\zeroh=x_0\lessdot x_1\lessdot \dots \lessdot 
x_n=\oneh$, 
and let $y$ and $z$ be two join-irreducibles with $\delta(y)=
\delta(z)$. Then $y$ and $z$ are incomparable.  
\end{lemma}

\begin{proof} Suppose on the contrary that $y<z$.  Let $j=\delta(y)=\delta(z)$.
Observe that $x_{j-1}\vee y =x_j \geq z$, so 
$$(z\wedge x_{j-1})\vee y = z\wedge(x_{j-1}\vee y)=z.$$
However, $z\wedge x_{j-1}$ and $y$ are both strictly less than $z$, so
$z$ is not join-irreducible, contrary to our assumption.  
\end{proof}

\section{Trim Lattices}

We now proceed to our study of trim lattices.  Let $L$ be a trim lattice,
with a specified left modular chain
$\zeroh=x_0\lessdot x_1\dots\lessdot x_n=
\oneh$.  

\begin{lemma} If $L$ is a trim lattice, 
it has exactly one join-irreducible and one 
meet-irreducible labelled $i$ for $1\leq i \leq n$. \end{lemma}

\begin{proof} Since $x_i$ is the join of the join-irreducibles labelled at 
most $i$, while $x_{i-1}$ is the join of the join-irreducibles labelled
at most $i-1$, there must be at least one join-irreducible labelled $i$.  
By trimness, there is exactly one.  The dual argument proves the statement
for meet-irreducibles.  \end{proof}

\begin{theorem} If $L$ is trim, so is any interval of $L$.
\end{theorem}

\begin{proof} If is sufficient to show that if $x\in L$, then
the interval $[\zeroh,x]$ is trim, since the dual result follows, and 
the trimness of $[y,x]$ can be proved by showing the trimness of $[\zeroh,x]$,
and then applying the dual result to the trim lattice $[\zeroh,x]$.

By Proposition 1, $[\zeroh,x]$ is left modular. 
Let the length of the 
left modular maximal chain in $[\zeroh,x]$ be $m$.  
We must show that there are
exactly $m$ join-irreducibles and $m$ meet-irreducibles in $[\zeroh,x]$. 

We consider $[\zeroh,x]$ labelled by the labelling induced from
$L$.
The join-irreducibles of $[\zeroh,x]$ are exactly the join-irreducibles
of $L$ that lie in $[\zeroh,x]$, and they have the same labels that they
do in $L$, so their labels are all different.  
Since, by Proposition 5, the labelling induced from $L$ agrees (up to an
order-preserving relabelling) with the labelling of 
$[\zeroh,x]$ induced from its left modular chain, 
the induced labelling uses $m$ different labels.  
Thus, 
$[\zeroh,x]$ has $m$ join-irreducibles, as desired.  

Let $a$ be a label that does not appear on a join-irreducible of 
$[\zeroh,x]$ (and which therefore doesn't appear in 
$[\zeroh,x]$ at all). 
Since the labelling on $[\zeroh,x]$ can also be considered as
being induced by its meet-irreducibles, there is no meet-irreducible of
$[\zeroh,x]$ labelled $a$.

Let $b$ be a label that appears on a join-irreducible of 
$[\zeroh,x]$.  Let $y$ be the join-irreducible of $L$ with that label
(which is also a join-irreducible of $[\zeroh,x]$).  Let $z$ be the
meet-irreducible of $L$ with label $b$.  
Let $\bar z=
z\wedge x$.  Since $\bar z \leq z$, and $y\not \leq z$, $y\not\leq \bar z$.
So $y\vee \bar z \ne \bar z$.  Let the increasing chain from 
$\bar z$ to $y\vee \bar z$ be $\bar z=t_0\lessdot t_1 \lessdot \dots
\lessdot t_r=\bar z\vee y$.  Since all the $t_i\leq x$, it follows that
$t_1 \not \leq z$ (otherwise $t_1\leq x\wedge z=t_0$, a contradiction).  
Thus, by the meet-irreducible labelling, $\gamma(t_0,t_1)\geq b$.  
By the join-irreducible labelling, $\gamma(t_{r-1},t_r)\leq b$.  Since
the labels on the chain are increasing, 
the chain consists of a single covering relation, which is labelled by $b$.

By the meet-irreducible labelling for $[\zeroh,x]$, it follows that
$\bar z$ lies below some meet-irreducible of $[\zeroh,x]$ labelled $b$.
But any element at the bottom of an edge labelled $b$ in $[\zeroh,x]$ lies
below $x$ and below $z$, thus below $\bar z$.  So $\bar z$ must be a
meet-irreducible labelled $b$ in $[\zeroh,x]$.  Since any other meet-irreducible
labelled $b$ in 
$[\zeroh,x]$ would have to lie below $\bar z$, $\bar z$ is the only one,
since two meet-irreducibles with the same label in a left modular lattice
must be incomparable, by the dual of Lemma 1.  Thus there is exactly one meet-irreducible 
labelled $b$, as desired.  \end{proof}

\begin{theorem}[{[Ma]}] If $L$ is trim and graded, it is distributive.
\end{theorem}

\begin{remark} A lattice with a chain of length $n$ (i.e. with
$n+1$ elements) and 
which has exactly $n$ join-irreducibles and $n$ meet-irreducibles is
called {\it extremal}.  Extremal lattices were introduced by
Markowsky in [Ma].  There, he showed that graded extremal lattices 
are distributive. Since trim lattices are by definition extremal, Theorem
2 follows. However, in the interests of self-containedness, we give a 
different proof.

It is worth noting that there are extremal lattices which are not trim.
Markowsky shows that any finite lattice can be embedded as an interval
of an extremal lattice, while Theorem 1 tells us that the intervals of
trim lattices are trim.  These two results imply that   
not all extremal
lattices are trim.  In particular, [Ma] gives an example of an extremal
lattice with 39 elements 
containing $M_3$ (see below)
as an interval; since $M_3$ is not trim, we know that 
this example is not trim.  
\end{remark}

\begin{proof} 
To show that a lattice is distributive, it suffices to show that
it has no sublattice $M_3$ or $N_5$ (see Figure 2) [Gr, Theorem II.1].  
This will follow
from the following two lemmas.

$$\epsfbox{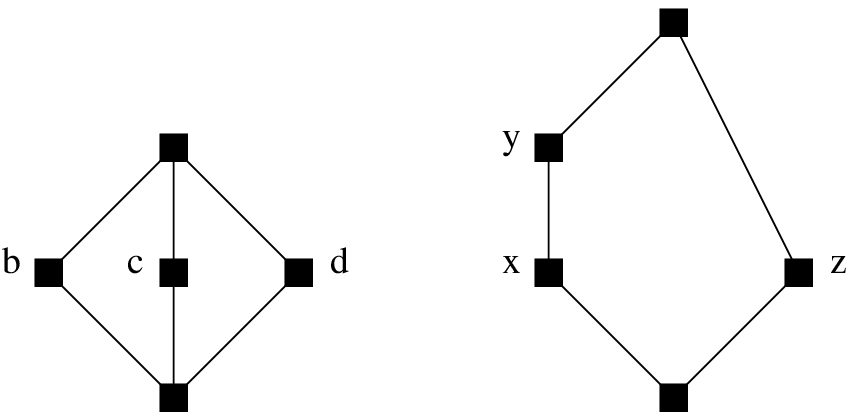}$$
$$\text{Figure 2: $M_3$ and $N_5$}$$
\begin{lemma} If $L$ is trim, then $L$ contains no 
sublattice isomorphic to $M_3$.
\end{lemma}

\begin{proof} We are going to assume that $L$ contains a sublattice
isomorphic to $M_3$ and obtain a contradiction.  
The property of being trim passes to intervals, by Theorem 1,
so we can reduce to the case where the minimum and maximum elements of
the copy of $M_3$ are $\zeroh$ and $\oneh$.  
Let the elements of the copy of $M_3$ be identified as in Figure 2.  

As always, 
let $n$ be the maximum label on the increasing chain from $\zeroh$ to $\oneh$.
Let $B$, $C$, $D$ denote the set of labels on the increasing chains from
$\zeroh$ to $b$, $c$, $d$ respectively.  
Suppose some two of them, say $B$ and $C$,
both contain $n$.  Then $b$ and $c$ both lie over some 
join-irreducible labelled $n$.  Since $b\wedge c=\zeroh$, there is more than
one join-irreducible labelled $n$, contradicting Lemma 2.  

On the other hand, suppose that some two of $B$, $C$, $D$, say $B$ and $C$, 
do not contain
$n$.  
Since $n$ is the maximum
label on the increasing chain from $\zeroh$ to $\oneh$, we can see from the
labelling induced from meet-irreducibles that every maximal chain from $\zeroh$ to
$\oneh$ contains an edge labelled $n$.  Thus, it occurs on both the increasing
chain from
$b$ to $\oneh$ and on the increasing chain from $c$ to $\oneh$.  We now apply
the dual of the previous argument to obtain a contradiction in this case
also.  

Since either two of $B$, $C$, $D$ contain $n$ or two do not, we are done.
\end{proof}

\begin{lemma} Let $L$ be a graded trim lattice.  Then $L$ contains no 
$N_5$.  \end{lemma}

\begin{proof} As in the proof of the previous lemma, 
we may assume that the minimum and maximum elements
of the $N_5$ are $\zeroh$ and $\oneh$.  
Let the other elements be 
identified as in Figure 2.  
Let $B$ be the set of labels on the increasing chain from
$x$ to $y$.  Suppose that the increasing chain from $\zeroh$ to $z$ has a 
label drawn from $B$, say $b$.  
Then $z$ and $y$ both lie over join-irreducibles 
with label $b$.  By the assumption that $L$ is trim, 
there is only one join-irreducible labelled $b$, 
so $z$ and $y$ both lie over it, which contradicts the assumption
that $z\wedge y=\zeroh$.  

Dually, no label from $B$ can occur on the increasing
chain from $z$ to $\oneh$.  
However, since we are assuming that $L$ is graded, the set of labels
appearing on every maximal chain from $\zeroh$ to $\oneh$ is the same, and
we have a contradiction.
\end{proof}

This completes the proof of Theorem 2.\end{proof}

\begin{theorem} If $L$ is trim, and $K$ is a sublattice of $L$ 
containing the left modular chain of $L$, then $K$ is trim.
\end{theorem}

\begin{proof} It is clear that the left modular chain in $L$ is still
left modular in $K$, so $K$ is left modular.  

Suppose $K$ is not trim.  Therefore, $K$ has either
two join-irreducibles with the same label, or two meet-irreducibles
with the same label.  Dualizing if necessary, we may assume that $K$
has two join-irreducibles with the same label, say $y$ and $z$, with 
label $b$.  
The fact that
$y$ and $z$ are labelled $b$ means that $y$ and $z$ lie below 
$x_b$ but not below $x_{b-1}$. This implies that, in $L$, each can 
be written as a join of join-irreducibles with labels no more than $b$,
and including at least one join-irreducible of $L$ labelled $b$.  Let
$j$ be the unique join-irreducible of $L$ labelled $b$.  So both
$y$ and $z$ lie over $j$.  It follows that $p=y\wedge z$ also lies
over $j$.  So $p$ lies below $x_b$ but not below $x_{b-1}$.  It follows
that, in $K$, $p$ lies over some join-irreducible labelled by $b$.  
But this implies that there are two 
join-irreducibles labelled by $b$ in $K$ which are comparable, and that
is impossible, by Lemma 1.  \end{proof}

\section{The sublattice fixed under a group of automorphisms}

The goal of this section is to show that if $L$ is a trim lattice, and a group
$G$ acts on
$L$ by lattice automorphisms, then $L^G$, the sublattice of $L$ consisting
of elements of $L$ fixed by $G$, is
a trim lattice.  

\begin{example} 
To orient oneself in this section, it is useful to consider the case where
$L$ is the Boolean lattice of all subsets of $[n]$, and $G=\{1,\sigma\}$
where $\sigma$ acts by interchanging $1$ and $n$.  

The first important thing to notice about this example 
is that the maximal chains in $L^G$ 
are shorter than the maximal chains in $L$.  The second thing to notice
is that if we make a reasonable-seeming choice of left modular chain by
setting $x_i=[i]$, only the top and bottom elements of our chosen
left modular chain actually lie in $L^G$.  Inspired by this example,
before we try to show that $L^G$ is trim, we will find some more 
left modular elements in $L$.  
\end{example}

For $L$ a trim lattice, we follow the terminology suggested by 
Drew Armstrong and 
say that the {\it spine} of $L$ 
consists of those elements of $L$ which lie on some chain of maximum length
in $L$.

\begin{lemma} If $L$ is a trim lattice, then all the elements of the
spine of $L$ are left modular.  
\end{lemma}

\begin{proof} This proof was suggested to me by Peter McNamara [Mc2].
Suppose $z$ is in the spine of $L$.  Let the labels which occur on the
increasing chain from $\zeroh$ to $z$ be $C$, and let the labels which
occur on the label from $z$ to $\oneh$ be $D$.  Since $z$ is 
in the spine, $C\cup D=[n]$, where $n$ is the length of the left modular
maximal chain in $L$.

Now suppose, for the sake of contradiction, that $z$ is not left modular.
It follows that there are some elements $p< q$ in $L$
such that $p\vee(z\wedge q)\ne
(p\vee z)\wedge q$.  Since the modular inequality tells us that 
$p\vee(z\wedge q)\leq (p\vee z)\wedge q$, it is in fact true that
$$p\vee(z\wedge q)<(p\vee z)\wedge q.$$
Now set $x=p\vee(z\wedge q)$, $y=(p\vee z)\wedge q$.  Note that 
$x\vee (z\wedge y)=x$, while $(x\vee z)\wedge y  =y$.  Thus, $x,y,z$ generate
a sublattice of $L$ of the following form:

$$\epsfbox{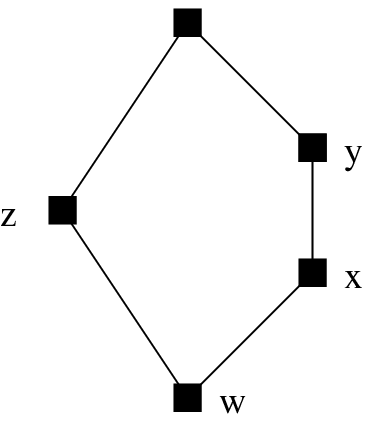}$$
$$\text{Figure 3: Sublattice of $L$ generated by $x,y,z$}$$

Let $b$ be a label on the increasing chain from $x$ to $y$.  Suppose
$b\in C$.  Since there is a unique join-irreducible labelled by
$b$, say $j$, it follows that $y$ and $z$ both lie over $j$.  But this
means that $w$ lies over $j$, and therefore $x$ lies over $j$, so $j$
cannot appear as a label on the increasing chain from $x$ to $y$, which 
contradicts our assumption.  
On the other hand, if $b\not \in C$, then $b\in D$, and
we can apply the dual argument to yield a contradiction.  
\end{proof}

The following lemma was suggested to me by Drew Armstrong [Ar], who observed
it to hold in the Tamari lattice.

\begin{lemma} The spine of a trim lattice $L$ is a distributive sublattice 
of $L$.
\end{lemma}

\begin{proof} The main difficulty is to show that the spine is closed under
lattice operations.  Suppose $y$ and $z$ are in the spine.  We will show that
$y\wedge z$ is also in the spine.  

Choose a left modular maximal chain $\zeroh=x_0\lessdot x_1\lessdot
\dots \lessdot x_n=\oneh$ in $L$.  
Let the set of
labels below $y$ be $A$, below $z$ be $B$, above $y$ be $C$, above $z$ be
$D$.  So $A\cup C=B\cup D=[n]$.  Let $P=A\cap B=\{p_1,\dots,p_r\}$,
with $p_1<\dots<p_r$.  Let $w_i$ be the join of the join-irreducibles 
indexed by $\{p_1,\dots,p_i\}$.  The $w_i$ are all distinct, since
$w_i$ lies below $x_{p_i}$ but not below $x_{p_i-1}$.  Thus, they form
a chain of length $r$ from $\zeroh$ to $y\wedge a$.  Since there are only
$r$ labels available for this chain (namely, the labels in $P$), and 
each label occurs at most once, this chain must be maximal.  

Now let $Q=C\cup D=\{q_1,\dots,q_{n-r}\}$, with $q_1<\dots<q_{n-r}$. 
Similarly to the above,
the meet-irreducibles lying above $y\wedge z$ are exactly those indexed by
$Q$.   
Let $v_{i}$ be the meet of the meet-irreducibles indexed by 
$\{q_i,\dots,q_{n-r}\}$.  By a dual argument, the $v_i$ are all distinct,
and therefore form a chain of length $n-r$ from $y\wedge z$ to $\oneh$,
which is necessarily maximal.  Thus $y\wedge z$ is in the spine of $L$,
and dually the same is true for $y\vee z$.  

We have shown that the spine is a sublattice of $L$.  By Lemma 3, it contains
no sublattice isomorphic to $M_3$, and by Lemma 5 it contains no 
sublattice isomorphic to $N_5$, so it is distributive.  
\end{proof}

We are now ready to prove the main theorem of this section:

\begin{theorem} If $L$ is a trim lattice and 
$G$ is a group which acts on $L$ by lattice automorphisms,
then $L^G$ is also trim. \end{theorem}

\begin{proof} First, we want to show that $L^G$ is left modular.  
Let $S$ be the spine of $L$.  Now the elements of $S^G$ are
left modular in $L^G$, since they are left modular in $L$.  To show
that $L^G$ is left modular, it remains to
show that $S^G$ contains a maximal chain in $L^G$.  It will suffice 
to show that if $y\lessdot z$ in $S^G$, then $y\lessdot z$ in $L^G$.  


Fix $y\lessdot z$ in $S^G$, and 
pick a maximal chain in $S$, 
$y=t_0\lessdot t_1\lessdot \dots \lessdot t_r=z$.
Let $j$ be the (unique) irreducible of 
$L$ which lies below $t_1$ but not below $y$.  Let the $G$-orbit of $j$ be
$\{j=j_1,\dots,j_k\}$.  Let $v_i=y\vee j_i$.   
Because $G$ acts by lattice
automorphisms, for every $i$, $y\lessdot v_i$, and $v_i\in S$.  
Let $w$ be the join of the
 $v_i$.  Observe that $w\in S^G$.  But $w\leq z$, so, since $y\lessdot z$ in
$S^G$, $z=w$.

The $v_i$ are all distinct, and since $S$ is distributive, the length of
any maximal chain in $S$ from $y$ to $z$ is of length $k$.  This means
that the only join-irreducibles lying below $z$ but not below $y$ are the
$j_i$.  

Now suppose that there is some $u$ in $L^G$ such that $y<u<z$.  There must
be some join-irreducible below $z$ but not below $y$ which is also below
$u$, but since $u\in L^G$, all the $j_i$ must lie below $u$, which would
force $u=z$, a contradiction.  

This implies that the maximal chains in $S^G$ are left modular maximal chains
in $L^G$ as desired.  

Now we want to show that $L^G$ is trim.  Let $T$ be the set of elements of 
$L^G$ formed by taking the join of the join-irreducibles in some 
$G$-orbit.  Clearly, any element of $L^G$ can be written as a join of
elements from $T$, so $T$ contains all the join-irreducibles of $L^G$.  
However, we showed above that if $y\lessdot z$ in $S^G$, then
there is exactly one $G$-orbit of irreducibles below $z$ but not below
$y$.  So there are chains in $L^G$ whose length is the number of 
$G$-orbits of irreducibles, which implies that there are at least that
number of join-irreducibles in $L^G$, so all the elements of $T$ are
join-irreducibles in $L^G$, and in particular, $L^G$ has the correct 
number of join-irreducibles to be trim.  Dually, $L^G$ has the correct 
number of meet-irreducibles, and it is therefore trim
\end{proof}

\section{The Level Condition and its consequences}

\begin{theorem} If $L$ is a trim lattice then it satisfies the
level condition of \emph{[BS]}:
$$\text{If $a$ and $b_1,\dots,b_k$ are atoms, and }
\delta(a)<\delta(b_1)<\dots<\delta(b_k),\text{ then } 
a\not<{b_1}\vee\dots\vee {b_k}.$$
\end{theorem}

\begin{proof} Suppose otherwise.  The proof is by induction on $k$.  
The statement is clearly true when $k=1$.  Suppose it is
true for $k-1$.  Consider a set of atoms $b_1,\dots,b_k$ as in the
statement of the theorem. 
Let $y={b_1}\vee\dots\vee {b_{k-1}}$,
and $z=b_2\vee \dots\vee b_k$.  
Since, by assumption, the statement is
true for $\{b_2,\dots,b_k\}$, ${b_1}\not\leq z$.
Since ${b_1},\dots,{b_{k-1}}$ all lie below $x_{\delta(b_{k-1})}$ while 
${b_k}$ does not, ${b_k}\not\leq y$.  

Suppose there is some atom $a$ with
$\delta(a)<\delta(b_1)$, such that $a<y\vee z$.  Choose such an $a$ with
$\delta(a)$ as small as possible.  Thus, we may assume that $a$ is the 
first element on the left modular chain from $\zeroh$ to $y\vee z$.  It follows
that $a$ appears on every maximal chain from $\zeroh$ to $y\vee z$.    
By the induction assumption,
$a$ lies below neither $y$ nor $z$.  Thus, $a$ appears as a label on 
the increasing chain from $y$ to $y\vee z$, and also on the increasing chain 
from $z$ to $y\vee z$.  Since the interval from $\zeroh$ to $y\vee z$ is
trim, there is some meet-irreducible in it labelled $a$, and both $y$
and $z$ lie below it.  But this contradicts the fact that $y\vee z$ is the
top of the interval.  Thus there can be no such atom $a$.\end{proof}

Recall that a lattice is said to be upper semimodular 
if, given three elements such
that $y$ and $z$ both cover $w$, then $y\vee z$ covers $y$ and $z$.  
Lower semimodularity is the dual condition.  A lattice is said to be
modular if it is both upper and lower semimodular.  Distributive
lattices are examples of modular lattices.  

Modularity implies gradedness, so we cannot hope that trim lattices will be
modular.  The following theorem shows that trim lattices posess a weakened
form of upper semi-modularity. 
The dual statement, which is also true, gives an
analogue of lower semi-modularity.    

\begin{theorem} Let $L$ be a trim lattice.  Let $y$ and $z$ cover 
$w$, and suppose that $\gamma(w\lessdot y)<\gamma(w\lessdot z)$. Then
$z\lessdot y\vee z$.  
\end{theorem}

\begin{proof} By Theorem 1, we can reduce to the case where 
$w=\zeroh$, $y$ and $z$
are atoms, and $y\vee z=\oneh$.  As usual, let $\zeroh=x_0\lessdot x_1
\lessdot \dots \lessdot x_r=\oneh$ be the left modular chain.  
By Theorem 5, $y$ is the 
join-irreducible
of $L$ with the smallest label, so $y=x_1$.  
The left modular maximal chain from $z$ to $\oneh$ consists of 
$z\vee x_i$.  But $x_1=y$, so the first element of this chain above $z$
is $y\vee z$.  Thus $y\vee z$ covers $z$.  
\end{proof}

We will call a lattice {\it nuclear} if $\oneh$ is the join of the atoms of 
$L$.  (In [Re], the term ``atomic interval'' is used for an interval 
in which the join of the atoms is the top of the interval.  Because this
might cause confusion with the standard use of atomic to describe a 
lattice in which every element can be written as a join of atoms,
we prefer to use a different term.)

\begin{theorem} If $L$ is trim and nuclear then its order complex
is homotopic to a sphere, whose dimension is 2 less than the number of
atoms of $L$.  
If $L$ is trim but not nuclear, then its order 
complex is contractible.  
\end{theorem}

\begin{remark} Note that since all intervals in a trim lattice 
$L$ are trim by Theorem 1, this theorem also applies to any interval in
a trim lattice.  \end{remark}

\begin{proof}
Observe that $x_1$ is the join-irreducible
labelled 1, and is an atom.  Any maximal chain in $L$ has an edge labelled
by 1; in a decreasing chain, this must be the last edge. The bottom of 
such an edge is a meet-irreducible labelled 1; thus, there is at most one
edge labelled 1 descending from $\oneh$.

Suppose $L$ is nuclear.  We prove the statement of the theorem by induction 
on the number of atoms of $L$.  
If $L$ has only one atom, the statement is obvious.  Suppose the 
statement holds for nuclear trim lattices with $r-1$ atoms.  
Let $a_1,\dots,a_r$ be the atoms of $L$, in 
increasing order by their labels.  Let $z=a_2\vee \dots \vee a_r$.  
The interval
$[\zeroh,z]$ is a nuclear trim lattice with $r-1$ atoms,
so by induction
it has a unique decreasing chain from $\zeroh$ to $z$.  This chain 
corresponds to a sphere of dimension $r-3$, so it is of length $r-1$.  
Now consider the 
increasing chain from $z$ to $\oneh=z\vee a_1$.  The top of this chain
is labelled with the label of $a_1$, which is 1, and all the other 
labels must be strictly greater than 1.  Since the chain is increasing,
this means that the chain is of length 1.  Thus, the decreasing chain
from $\zeroh$ to $z$ extends uniquely to a decreasing chain from $\zeroh$ to $\oneh$.
By the remarks in the first paragraph, any decreasing chain from $\zeroh$ to 
$\oneh$ 
passes through
$z$.  Since there is only one decreasing chain from $\zeroh$ to $z$, the 
decreasing chain form $\zeroh$ to $\oneh$ which we have found is unique,
and it is clearly of length $r$, which implies that the order complex of 
$L$ is homotopic to a sphere of dimension $r-2$, as desired.    

For the second statement, it
 is well-known that if $L$ is any finite non-nuclear lattice, then
its order complex is contractible.  This follows from the Crosscut Theorem;
see, for example [B3].  
\end{proof} 

One of the reasons to be interested in statements about the homotopy types
of order complexes of intervals is that for $x<y$ in any poset,
the M\"obius function $\mu(x,y)$ is the reduced Euler characteristic of
the order complex of the interval $[x,y]$.  Thus, from Theorem 7
combined with Theorem 1, we
can deduce that the M\"obius function of any interval in a trim lattice
is either 0, 1, or $-1$.  
That $\mu(\zeroh,\oneh)$ is 0, 1 or $-1$ for a trim lattice also follows easily 
from results in [BS].    

\section{Cambrian Lattices}

Let $W$ be a finite subgroup of the orthogonal transformations of 
$E=\Bbb{R}^n$, generated by reflections.  Such a group is called a 
{\it finite reflection group}.
It has an associated finite root system
$\Delta\subset E$, 
which is partitioned into positive and negative roots, denoted
$\Delta^+$ and $\Delta^-$.  The elements of $W$ permute $\Delta$.  

For any element $w$ of $W$, let the {\it inversion set} of $w$ be defined by:
$$I(w)=\{\alpha\in \Delta^+\mid  w^{-1}(\alpha)\in \Delta^-\}.$$

If we order the elements of $W$ by inclusion of inversion sets, we obtain
a partially ordered set structure 
called {\it weak order} on $W$.  Weak order on $W$ is a 
lattice.  A general reference for weak order on finite reflection 
groups is [B2].

A lattice homomorphism is a map of lattices which preserves lattice 
operations.  A quotient lattice of a lattice $L$ is the image of a 
homomorphism from $L$.  
The fibers of a lattice homomorphism from $L$ are 
necessarily intervals in $L$.  

Associated to any finite reflection group $W$ is a graph called its Coxeter
diagram, which we denote $G$.  
Let $\bar G$ be an orientation of $G$ (that is to say,
for each edge of $G$, we designate one end of the edge as the source and
the other as the target).  Associated to $\bar G$ is a {\it Cambrian lattice}
$C(\bar G)$, which is a quotient of weak order on $W$.  
We shall not give the general definition here, restricting our attention to
reflection groups of types $A$ and $B$, where (in contrast to the other types) 
explicit descriptions of
the Cambrian lattices are known.  
The general definition and the explicit description in types $A$ and $B$ are
due to Reading [Re].  In what follows, we will review these descriptions,
and then show that the Cambrian lattices in types $A$ and $B$ are trim,
and consequently that the results of the first half of this paper apply
to them (and their intervals).  The result of Theorem 7 applied to
Cambrian lattices was already proved in [Re].  Theorems 4 and 5 are new.

\subsection*{Type $A$ Cambrian Lattices}
In type $A_{n-1}$, the reflection group $W$ is isomorphic to $S_n$. 
Let $e_1,\dots,e_n$ be a basis for $\Bbb{R}^n$.  A permutation $\pi \in S_n$ 
acts on
$\Bbb{R}^n$ by taking $e_i$ to $e_{\pi(i)}$.  The roots are the vectors
$e_j-e_i$ for $i\ne j$.  The positive roots are $e_j-e_i$ for $j>i$.  

For $\pi \in S_n$,  
$e_j-e_i$ is an inversion of $\pi$ for $j>i$ if $j$ precedes $i$ in the word
$\pi_1,\dots,\pi_n$. As already mentioned, weak order on $S_n$ is 
the inclusion order on inversion sets.

The Coxeter diagram $G$ consists of a path of $n-1$ nodes, labelled 
$s_1, \dots, s_{n-1}$.  Let $\bar G$ be an orientation of this diagram.
We write $s_{i-1}\rightarrow s_i$ and $s_{i-1}\leftarrow s_i$ to 
represent the two possible orientations of the edge between $s_{i-1}$ and
$s_i$.  
Define two complementary subsets of $[2,n-1]$ by
$D = \{i\mid s_{i-1}\rightarrow s_i\}$, $U=\{i\mid s_{i-1}\leftarrow
s_i\}$. 

For our purposes, a {\it pattern} is a permutation of $[k]$. A permutation
$\pi$ 
contains a given pattern $\sigma$ if there are some $i_1<\dots<i_k$ such
that $\pi_{i_1},\pi_{i_2},\dots,\pi_{i_k}$ are in the same relative
order as $\sigma_1,\dots,\sigma_k$.  If we put a bar over an element of a 
pattern (as, for example, in $\bar 231$), then to say that $\pi$ contains
that pattern means that $\pi$ contains an instance of the pattern in 
which the element of $\pi$ that corresponds to the barred element of the
pattern
belongs to $U$.  
Similarly,
if we underline an element of the pattern, we mean that the corresponding
element of $\pi$ must belong to $D$.  

Let $\mB$ be the set of permutations in $S_n$ avoiding $\bar 231$ and
$31\underline2$.  
Let $\mT$ be the set of permutations in $S_n$ avoiding
$\bar 213$ and $13\underline2$.

\begin{proposition}[{[Re]}] There is a quotient of weak order on
$S_n$ 
the minimal elements of whose fibers are $\mB$ and the maximal elements of
whose fibers are $\mT$.  This is by definition the Cambrian lattice $C(\bar G)$. 
$\mB$ and $\mT$ are sublattices of weak order on $S_n$, each also isomorphic to
$C(\bar G)$.  
\end{proposition}

Note that in the case that all the edges of $\bar G$ are oriented 
$s_{i-1}\rightarrow s_i$, 
$\mB$ consists
of all those
permutations avoiding $312$, while $\mT$ consists of those permutations avoiding
$132$.  In this case, $C(\bar G)$ is a Tamari lattice, and the map from
$S_n$ to $C(\bar G)$ is the well-known quotient map from weak order on 
$S_n$ to the Tamari lattice.  (See, for instance, [BW] for more details.)

Here is an example, showing weak order on $S_3$, an oriented Dynkin diagram,
and the induced Cambrian lattice.  
$$\epsfbox{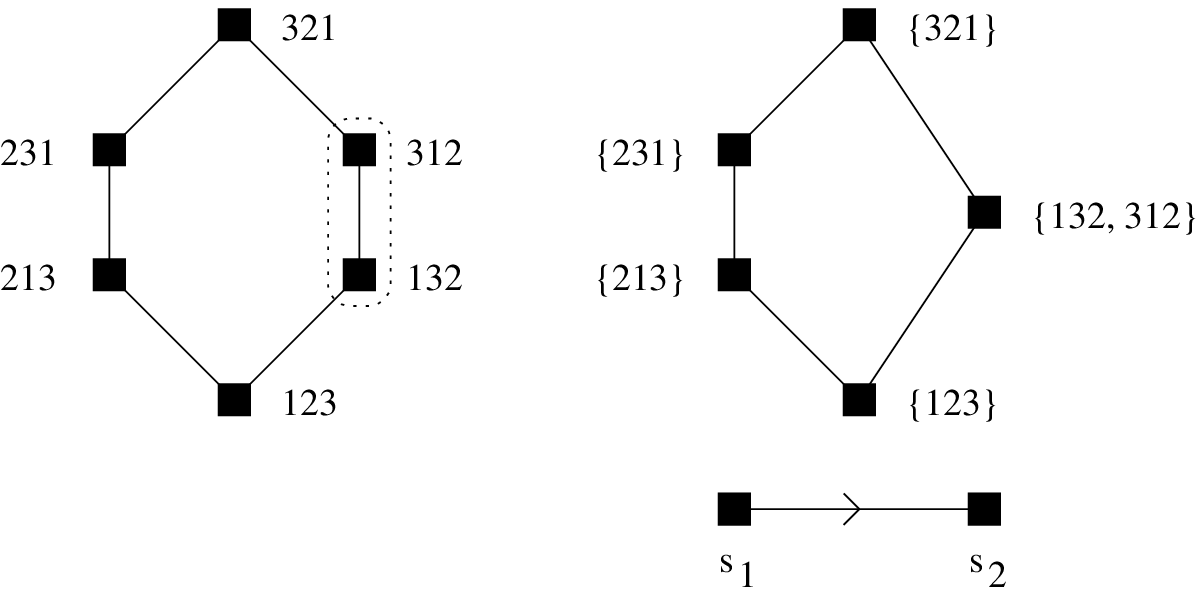}$$
$$\text{ Figure 4}$$

\subsection*{Type $B$ Cambrian lattices} 
We now consider the type $B_n$ Cambrian lattices.  Here, $W$ is isomorphic 
to $B_n$, the 
group of signed permutations of $[n]$, that is, permutations of 
$\{-n,\dots,-1,1,\dots,n\}$ which are fixed under the involution 
interchanging positive and negative numbers.  We think of $\pi \in B_n$ 
as acting
on $\Bbb{R}^n$ by taking $e_i$ to $e_{\pi(i)}$, where we let $e_{-i}=-e_i$.
The roots of $B_n$ are $\pm e_i \pm e_j$ for $i\ne j$, together with
$\pm e_i$.  The positive roots are those of the form 
$e_j-e_i$ for $j>i$, $e_j+e_i$, and $e_i$.  

The Coxeter diagram $G$
consists of a path of $n$ vertices, labelled $s_0,\dots,s_{n-1}$,
where the edge from $s_0$ to $s_1$ is the special edge labelled 4.  
(For those unfamiliar with Coxeter diagrams, this encodes certain information
about the relations among the 
generators of the group which correspond to these nodes, but this
is not essential for our purposes.)
Let $\bar G$ be
an orientation for $G$.  We define two complementary subsets of 
$\{-n+1,\dots,-1,1,\dots, n-1\}$: for $1\leq i\leq n-1$, 
if $s_{i-1}\rightarrow s_i$, then $i\in D$ and $-i\in U$, and vice versa
if $s_{i-1}\leftarrow s_i$.  

The one-line notation for $\pi \in B_n$ is 
$\pi_{-n}\dots\pi_{-1}\pi_1\dots\pi_n$.
For $0<i<j$, $e_j-e_i$ is an inversion of $\pi$
if $j$ precedes $i$ in the one-line notation for $\pi$,  
$e_i$ is an inversion of $\pi$ if $i$ precedes $-i$ in
the one-line notation for $\pi$, and 
$e_j+e_i$ is an inversion for $\pi$ if $i$ precedes $-j$ (or equivalently
$j$ precedes $-i$)
in the one-line notation for $\pi$.  
Weak order on $B_n$ is defined (as always) by
inclusion of inversion sets.  

As in type $A$, we say that $\pi$ contains a pattern $\sigma \in S_k$
iff there are some $i_1<\dots<i_k$ such that the relative order of 
$\pi_{i_1},\pi_{i_2},\dots
\pi_{i_k}$ is the same as that of $\sigma_1,\dots,\sigma_k$ --- 
but we allow $i_1,\dots,i_k$ to be chosen
from $[n]\cup -[n]$.  The meaning of overlines and underlines are the same
as in type $A$.  (Note: sometimes $B_n$ is considered as a set of
permutations on $\{1,2,\dots,n,\bar 1,\bar 2,\dots,\bar n\}$, i.e.\ $
\bar k$ is used where we would write $-k$.  There is a possibility of
confusion for the reader here.  When we write $\bar k$ we will never mean
$-k$; we will always mean that $k\in U$.)

The definition of the Cambrian lattice $C(\bar G)$ is very similar to the 
definition in type $A$.    
Let $\mB$ be the set of permutations in $B_n$ avoiding $\bar 231$ and
$31\underline2$.  (Note that, because of the symmetry of elements of $B_n$, it 
is actually sufficient to check that an element of $B_n$ avoids one of 
these patterns --- the avoidance of the other pattern comes for free.)
Let $\mT$ be the set of permutations in $B_n$ avoiding
$\bar 213$ and $13\underline2$.  (As for $\mB$, we only need to check one of 
these conditions.)  An analogue of Proposition 6 holds in type $B_n$:

\begin{proposition}[{[Re]}] There is a quotient of weak order on
$B_n$ 
the minimal elements of whose fibers are $\mB$ and the maximal elements of
whose fibers are $\mT$.  This is by definition the Cambrian lattice $C(\bar G)$. 
$\mB$ and $\mT$ are sublattices of weak order on $B_n$, each also isomorphic to
$C(\bar G)$.  
\end{proposition}

\subsection*{ Trimness of Cambrian lattices}
The story in type $B$ is in some respects
simpler than in type $A$, so we begin with the
following theorem:

\begin{theorem} 
The type $B$ Cambrian lattices are trim.\end{theorem}

\begin{proof}
First, we must understand the join-irreducibles of $C(\bar G)$.  

\begin{lemma} There are $n^2$ join- and meet-irreducibles of $C(\bar G)$.
\end{lemma}

\begin{proof} 
A 
join-irreducible $\pi$ of $\mB$ is necessarily a 
join-irreducible of $B_n$, because $C(\bar G)$ is a quotient of $B_n$.
Let the unique element which lies immediately 
below $\pi$ in $B_n$ be $\sigma$.  
Let the adjacent transposition relating $\pi$ and $\sigma$ interchange 
$x$ and $y$, (and also $-x$ and $-y$), with $y>0$, and $|y|\geq |x|$.  
(So $x$ and $y$ appear
together, in the order $xy$ in $\sigma$, and in the order $yx$ in $\pi$.)
We wish to show that $x$ and $y$ determine $\pi$.  

We consider first the case where $0<x<y$.  
Thus $\pi$ looks like  either 
$$\dots (-x)(-y) \dots yx \dots \text{ or } \dots yx \dots (-x)(-y) \dots$$
The fact that $\pi$ is join-irreducible in $B_n$ means that each of the three
segments into which $\pi$ is divided by $yx$ and $(-x)(-y)$ must be
increasing.  This immediately rules out the second of the two possibilities
displayed above.
Again using the fact that each of the three segments of $\pi$ is increasing,
to show that $\pi$ is determined by $x$ and $y$, it suffices to show that
for any $z$ other than $x$, $y$, $-x$, or $-y$, we can determine which 
segment it belongs to.  If $z>y$, then $z$ must occur in the rightmost 
segment.
If $0<z<x$, then $z$ must lie in the middle segment.  If $x<z<y$, then
$z$ cannot lie in the leftmost segment, and which of the other two segments
it lies in is determined by the fact that $\pi\in \mB$, and thus that 
exactly one of 
$zyx$ or $yxz$ is a forbidden configuration (depending on whether 
$z\in U$ or $z\in D$).   This determines the position of all $z>0$, and 
by the symmetry of $\pi$, it also determines the positions of all $z<0$.
Thus we see that $\pi$ is determined by $x$ and $y$.  

The cases where $x=-y$ and where $-y<x<0<y$ are very similar.  Thus, for
every pair $x,y$ with $-y \leq x <y$, there is exactly one join-irreducible
in $\mB$, and thus there are $n^2$ in total.  
Using a dual argument, there are exactly $n^2$ 
meet-irreducibles of $\mT$. Using Proposition 7 which says that 
$\mB$ and $\mT$ are both isomorphic to $C(\bar G)$, we see that 
$C(\bar G)$ has exactly $n^2$ join-
and meet-irreducibles.\end{proof}

Our next step will be to identify a maximal chain of length $n^2$ in 
$C(\bar G)$.  

We let $s_i$ (the labels of the nodes of the Coxeter diagram) also 
denote the corresponding reflection in $B_n$: for $i>0$, $s_i$ interchanges
$i$ and $i+1$, while $s_0$ interchanges $1$ and $-1$.  

Write out a word in which each $s_i$ occurs once, and such that for any edge
$s_i\rightarrow s_j$, $s_i$ occurs to the right of $s_j$ in the word.  
Let $c$ be the 
product of the $s_i$ in this order.  
It is a Coxeter element, and one
convinces onself easily that it takes $-n$ to the smallest element of $D$,
each element of $D$ to the next largest one, the largest element of $D$ 
to $n$, and by symmetry $n$ to the largest element of $U$, etc.  

We know that $c^n = -1$.  (This holds for any
Coxeter element in type $B$.  It is also easy to see
from our explicit description.)  
For $0\leq i \leq n^2$, let $x_i$ denote the element of $B_n$
which consists of the product of the rightmost $i$ 
simple reflections in $c^n$ (where we 
think of $c$ as being written as a word of length $n$ as above).  
Since the minimum length of an expression for $-1$ as a product of 
simple reflections in $B_n$ has length $n^2$,
and our expression for $c^n$ has exactly this length, 
it follows that $x_i\lessdot x_{i+1}$ in $B_n$.  

\begin{lemma} The $x_i$ are contained in $\mT\cap \mB$ (so in particular,
each determines a different element of $C(\bar G)$).  \end{lemma}

\begin{proof} In order to prove this, we will need to give a description
of inversion sets of elements of $\mT$ and $\mB$ in terms of their allowed
intersections with irreducible 
rank 2 root systems contained in our $B_n$ root system $\Delta$.

These are the types of rank 2 root systems contained in $\Delta$:

(i) The type $B_2$ root system corresponding to positive
roots $e_i, e_j+e_i, e_j, e_j-e_i$,
$(i<j)$.  To read off which of these elements lie in the inversion set of 
$\pi\in B_n$, we need only consider the relative positions of 
$i,j, -i, -j$.

(ii) The type $A_2$ root system corresponding to positive roots
$e_j-e_i, e_k-e_i,e_k-e_j$ for $i<j<k$.  To read off which of these elements lie
in the inversion set of $\pi$, we need to look at the relative positions of
$i,j,k$.   

(iii)  The type $A_2$ root system corresponding to positive roots
$e_j-e_i, e_j+e_k, e_i+e_k$ (for $i<j$).    To read off which of these elements
lie in the inversion set of $\pi$, we need to look at the relative positions of
$i,j,-k$.

The inversion set for any element of
$B_n$ intersected with any of these rank 2 root systems is an initial or
final subset of the list of inversions (in the order in which they are
listed above).  This can be seen by inspection in our case; a similar
statement holds for all finite reflection groups, see [B2].

\begin{lemma} For $R$ a rank 2 root system contained in $\Delta$, there 
is an order on its roots (either the one given above or its reverse) which
we call the $\bar G$-order such that:

(i) the inversion set of an element of 
$\mB$ intersected with $R$ is either an initial
subset with respect to the order, or consists of exactly the final element.

(ii) the inversion set of an element of $\mT$ intersected 
with $R$ is either an initial
subset or consists of all the elements except the first.  
\end{lemma}

\begin{proof} This essentially follows by inspection, considering the three
possible types of root systems contained in $B_n$.  
Suppose that the rank 2
root system is of type $B_2$.  The possible relative positions for $i,j,-i,
-j$ (ignoring all other symbols) are as follows:
$$\begin{array}{rlccrl} \multicolumn{3}{r}
{ji(-i)(-j),}& \multicolumn{3}{l}{I=\{e_i, e_j+e_i, e_j, e_j-e_i\}} \\
ij(- j)(-i),& I=\{e_i, e_j+e_i, e_j\}&  & &
j(-i) i (-j),& I=\{e_j+e_i,e_j, e_j-e_i\}\\
i(-j) j (-i), & I=\{e_i, e_j+e_i\} & & &
 (-i) j (-j) i,& I=\{e_j, e_j-e_i\} \\
(-j) i (-i) j, & I=\{e_i\} & & &
(-i) (-j) j i,& I=\{e_j-e_i\} \\
\multicolumn{3}{r}{(-j) (-i) i j,}&\multicolumn{3}{l} 
{I=\emptyset} \end{array}$$

Observe that if $i\in D$ then $ij(-j) (-i)$ and $i(-j) j (-i)$
are impossible for an element of $\mB$, 
while if $i\in U$ then $j(-i) i (-j)$ and 
$(-i) j (-j) i$ are impossible for an element of $\mB$.  Thus, if $i\in D$, 
part (i) of the 
lemma is satisfied if we set the $\bar G$ order to be $e_j-e_i, e_j,
e_j+e_i, e_i$, while if $i\in U$, part (i) of the 
lemma is satisfied if we set
$\bar G$-order to be the reverse order.  It is straightforward to check that
the same order also satisfies part (ii) of the lemma.  

The other two types of root systems are handled similarly, proving the 
lemma.
\end{proof}

We now prove a converse to Lemma 9.  First, we introduce some notation.
We say that a subset of a rank 2 root system is {\it initial} if it is 
initial with respect to the $\bar G$-order.  
We say that the subset is {\it last} if 
it consists of only the final element (with respect to the $\bar G$-order).  
We say that a subset is {\it all but first} if it consists of all the elements
except the first.  
We will say that a set of roots has $\mB$-good intersection with a 
rank two
root system if its intersection is initial or last, and 
$\mT$-good intersection if its intersection is initial or all but first.  
Thus, Lemma 9 says that if $w\in \mB$ then $I(w)$ has $\mB$-good intersection
with every rank 2 root system in $\Delta$, and similarly with $\mT$ replacing
$\mB$. The following lemma is a converse.  

\begin{lemma} If a set of roots has 
$\mB$-good intersection with every rank 2 root system, 
then the set of roots
is the inversion set of an element of $\mB$.  Similarly, if it has
$\mT$-good intersection with every rank 2 root system, it is the inversion
set of an element of $\mT$.  
\end{lemma}

\begin{proof} We prove the first statement.  
Given a set of roots $I$ whose intersection with any rank 2 root
system is either initial or final, it is the inversion set of a unique element
$\pi$ 
of $B_n$ [B2].  We must show that $\pi$ 
contains neither a $\bar 231$ nor a $31\underline2$. Suppose it does,
and  suppose first that 
this pattern involves three elements of distinct absolute values.  
If these are all the same sign (which we may assume to be positive) then
we have found $i,j,k$ such that $e_i-e_j, e_i-e_k,e_j-e_k$ has 
an illegal intersection with $I$.  Similarly, if the pattern involves 
elements not all of the same sign then we are in the other type of $A_2$ root
system, while if the pattern involves two elements of the same absolute 
value, 
then we are in a similar situation with respect to a $B_2$ root system.  
Thus, $\pi$ contains no $\bar 231$ or $31\underline2$, and therefore is
an element of $\mB$.  

The second statement follows from a similar argument, and the lemma
is proved.  \end{proof}

Now that we understand the possible inversion sets of element of $\mB$ and
$\mT$, we can return to the proof of Lemma 8.  
We now proceed to show that the inversion set of $x_i$ intersected with any
rank two root system is initial.  Consider, for example, a root system of
type $B_2$.  We must determine whether, in our word for $c^n$, the inversions
$e_i, e_i+e_j, e_j, e_j-e_i$ appear in that order or the reverse order.
By inspection (recalling our explicit description of $c$), 
we see that they occur in the forward order if $i\in U$, and
in the backward order if $i\in D$.  But now observe that this order on 
the roots is exactly the order provided by Lemma 9, as desired.  The other
types of root systems are dealt with similarly.  This completes the proof
of Lemma 8.  \end{proof}




We wish to show that the $x_i$ are left modular.  
By Proposition 4, it is sufficient to exhibit an interpolating labelling
for $C(\bar G)$ such that the $x_i$ form the increasing chain from
$\zeroh$ to $\oneh$.  

We now introduce some notation related to $C(\bar G)$.  For $x\in B_n$,
we write $[x]$ for the fibre of the quotient map to $C(\bar G)$ which
includes $x$.  We also write  $\piu(x)$ for the top element of $[x]$,
and $\pid(x)$ for the bottom element of $[x]$.  

We define an edge-labelling for $C(\bar G)$ as follows.  First,
observe that the edges in the Hasse diagram of weak order on $B_n$ have
a natural labelling by positive roots: we label the edge $x\lessdot y$
by $I(y)\setminus I(x)$.  We now use this labelling to define a labelling
for $C(\bar G)$.

Suppose 
$[x]\lessdot [y]$ in $C(\bar G)$, with $x\in \mT$.  Then $x$ is covered by
an element of $[y]$, say $y'$.  Then set  
$\gamma([x]\lessdot [y])=I(y')\setminus I(x)$.

\begin{lemma} If $x\lessdot y$ and $[x]\ne [y]$, then
$\gamma([x]\lessdot [y])=I(y)\setminus I(x)$.  
\end{lemma}

\begin{proof} Let $x'=\piu(x)$, and let $y'$ be the element of $[y]$ covering
$x'$.   
Now $y'\wedge y$ is in $[y]$ but lies over $x$, so must equal $y$, which
implies that $y\leq y'$.  
Since
$y$ does not lie under $x'$, but $x$ does, $I(y)\setminus
I(x)=I(y')\setminus I(x')$, as desired.  \end{proof}

We now prove an easy lemma which will be useful for computations in 
$C(\bar G)$.  

\begin{lemma} 
If $x,y \in \mT$, then $I(x\vee y)=I(x)\cup I(y)$.  
If $x,y\in \mB$, then $I(x\wedge y)=I(x)\cap I(y)$.  
\end{lemma}

\begin{proof} We prove the first statement.  Observe that $I(x)\cup I(y)$ 
has $\mT$-good intersection with every rank 2 root system, and therefore,
by Lemma 10, defines an element of $\mT$.  Now it is clear that this element
must be the join of $x$ and $y$. 

The argument for the second statement is similar. 
\end{proof}

\begin{lemma} The labelling $\gamma$ defined above is an
interpolating labelling for $C(\bar G)$.  \end{lemma}

\begin{proof}
The first necessity for showing that a labelling is interpolating 
is to show that it is an EL-labelling.  Let $[v]<[w]$ in $C(\bar G)$,
with $v$ and $w$ in $\mT$.  
Let
$\alpha=\min(I(w)\setminus I(v))$.

We begin by showing that there is a $z$ such that $[v]\lessdot [z] \leq [w]$, 
with
$\gamma(v\lessdot z)=\alpha$. 
Let $x$ be the element of $\mB\cap \mT$ whose inversion set consists of all roots
up to and including $\alpha$ in $\bar G$-order.  
Let $z=v\vee x$.  By Lemma 26, $I(z)=I(v)\cup I(x)=I(v)\cap\{\alpha\}$.
Thus $\gamma([v]\lessdot [z])=\alpha$, and clearly $[v]\lessdot [z]$ is
the only edge proceeding up from $[v]$ labelled by $\alpha$.  

Next we show that every maximal chain from $[v]$ to $[w]$ has an edge
labelled by $\alpha$.  Given a maximal chain, let 
$[q]$ be the first element of the chain lying over $z$, 
and let $[r]$ be the element lying below $[q]$ in the chain.
Let $r\in \mT$.  Let $q'$ be the element of $[q]$ covering $r$.  Then 
$q'$ lies over $z$ but $r$ does not, so $I(q')\setminus I(r)=
\{\alpha\}$, and the edge $[r]\lessdot [q]$ is labelled $\alpha$.  

So $\alpha$ is the minimum possible label to occur on any edge of any
maximal chain from $[v]$ to $[w]$, and it must occur on every chain.  
Thus, the first step in any increasing chain from $[v]$ to
$[w]$ must be labelled $\alpha$, so any increasing chain must begin
$[v]\lessdot[z]$.  Now, by induction, there is a unique increasing chain
from $[v]$ to $[w]$.  

Now we must show that the labelling $\gamma$ is interpolating.  
So suppose that we have chain of length two which isn't increasing,
say $[v]\lessdot[u]\lessdot[w]$.  Let us assume that $v \in \mT$.  
Let $\alpha=\gamma([v]\lessdot[u])$ and $\beta=\gamma([u]\lessdot[w])$.  
Let $[v]=[y_0]\lessdot [y_1]\lessdot \dots\lessdot [y_r]=[w]$ be the 
increasing chain from $[v]$ to $[w]$.  

Since $\gamma([y_0]\lessdot[y_1])$ is the minimum label on any chain in the
interval, by what we have just shown it must occur on every chain from
$[v]$ to $[w]$.  It cannot be that the edge $[v]\lessdot [u]$ has this
label, so $\gamma([y_0]\lessdot[y_1])=\gamma([u]\lessdot [v])$, one of
the two conditions necessary for $\gamma$ to be interpolating.  

In weak order on $B_n$, we know that there are two edges rising from 
$v$, labelled by $\alpha$ and $\beta$.  These correspond to simple
reflections $s_\alpha$ and $s_\beta$ (i.e. the tops of these edges 
are $vs_\alpha$ and $vs_\beta$ where $s_\alpha$ and $s_\beta$ are simple
reflections.) Let $V$ be the subgroup of $W$ generated by 
$s_\alpha$ and $s_\beta$.  Then $v$ is the unique minimum-length 
representative of its left coset $vV$ in $W$.  This coset appears in
weak order on $B_n$ as an interval with minimum element $v$.  (For 
more details, see [Hu, Section 1.10].)  This interval of $B_n$ is
isomorphic to weak order on $V$, which is a rank 2 reflection group.  
Therefore, this interval consists of a two incomparable chains $
vs_\alpha=c_1
\lessdot c_2\lessdot \dots \lessdot c_k, vs_\beta=d_1\lessdot d_2\lessdot
\dots \lessdot d_k$, 
together
with a minimum element $v$ and a maximum element, which we will call
$q$.  

Observe that $c_1=vs_\alpha\in [y_1]$ and $d_1=vs_\beta\in [u]$.
Thus, their
join, $vs_\alpha\vee vs_\beta$, which equals $q$, lies in $[w]$.  
Since $[q]>[d_1]$, but $c_k \not >d_1$, $[c_k]\ne [q]$.
Since $[y_1]=[c_1]\leq [c_k]\lessdot[q]=[w]$, $[c_k]=[y_{r-1}]$.  

Observe that the edge (in weak order on $B_n$) from $c_k$ to $q$ is labelled
by $\beta$.  Thus, by Lemma 11, since $c_k\in [y_{r-1}]$ and
$q\in [w]$, $\gamma([y_{r-1}]\lessdot[w])=\beta$, and
we have shown that $\gamma$ is interpolating.  
\end{proof}

It is clear that the $x_i$ form the 
increasing chain from $\zeroh$ to $\oneh$ in $C(\bar G)$, and thus
they are left modular.  
We conclude that 
Cambrian lattices of type $B$ are trim.  \end{proof}

The following theorem is an easy corollary of Theorem 8.

\begin{theorem} The type $A$ Cambrian lattices are trim.   \end{theorem}

\begin{proof}
Let $\bar G$ be an oriented type $A$ Coxeter diagram.  Let $\bar G'$ be 
the type $B$ Coxeter diagram obtained by affixing an extra edge 
labelled 4 to 
$G$, oriented arbitrarily.  It is straightforward to see, either by 
the explicit description of Cambrian lattices in types $A$ and $B$, or 
from general theory, that $C(\bar G)$ is a lower interval in
$C(\bar G')$.  (The top of the interval is the equivalence class of the
longest word for the type $A$ Coxeter group.)  It now follows by Theorem 1 
that 
$C(\bar G)$ is trim.  
\end{proof}

\subsection*{Conjectural description of other Cambrian lattices}
Let $W$ be a finite reflection group which contains $-1$.  Let $G$ be
its Dynkin diagram, and $\bar G$ an orientation.  As in type $B$, we 
can order the nodes of the diagram in accordance with the orientation
of the edges, and then take the product of the simple reflections
in that order, obtaining 
a Coxeter element $c$.  If $h$ is the Coxeter number for $W$ then,
since $-1 \in W$, $h$ will be even, and $c^{h/2}=-1$, [Hu, Corollary 3.20].

Linearly order the roots of $W$ in the order in which they appear as 
inversions in the word for $c^{h/2}$.  Let $x_i$ be the element of 
$W$ whose inversion set consists of the first $i$ roots, in this order.  

Now, take the minimal quotient of $W$ such that the $x_i$ are all 
left modular.  Call this the pre-Cambrian lattice associated to $\bar G$.  

\begin{conjecture} The bottom elements of the fibres of this quotient will
be exactly those elements whose inversion sets have $\mB$-good intersection
with all rank 2 sub-root systems, where the order on the sub-root system
comes from the linear order on the positive roots.  (And similarly for 
the top elements of the fibres.) \end{conjecture}

\begin{conjecture} The pre-Cambrian lattice associated to $\bar G$
coincides with the Cambrian lattice $C(\bar G)$.  \end{conjecture}

Note that we 
have already showed that these conjectures hold in type $B$.

Since every root system embeds in one whose reflection group contains
$-1$, and the Cambrian lattice associated to the smaller root system
appears as a lower interval in the Cambrian lattice associated to the larger
root system, it would follow from Conjectures 2 and 3 
that all Cambrian lattices
are trim.  

\section*{ Acknowledgements}

I would like to thank Peter McNamara for discussions in which the fundamental
ideas of this paper were formed, and for many helpful comments as this project
developed.  I would also like to thank Drew Armstrong,  
Nathan Reading, Vic Reiner, 
Bruce Sagan, and an anonymous referee for 
their useful comments and suggestions.

\end{document}